\documentclass[12pt]{amsart}
\usepackage{graphicx}
\usepackage{amssymb}
\usepackage[shortlabels]{enumitem}
\usepackage[margin=1.25in]{geometry}
\newtheorem{thm}{Theorem}[section]

\newtheorem{prop}[thm]{Proposition}

\theoremstyle{definition}

\theoremstyle{remark}
\newtheorem{rem}[thm]{Remark}

\begin{document}
\title[]{Hilbert's $19^{\text{th}}$ Problem Revisited}
\author[]{Connor Mooney}
\maketitle
\tableofcontents

\newpage
\section{Introduction and Acknowledgments}
In this survey article we revisit Hilbert's $19^{\text{th}}$ problem concerning the regularity of minimizers of variational integrals. Sections $2-4$ are devoted to the classical theory (that is, the statement and resolution of Hilbert's problem in all dimensions). In sections $5-6$ we discuss recent results concerning the regularity of minimizers of degenerate convex functionals. In the last section we discuss some open problems. Exercises are included for the benefit of researchers who are entering the subject.

\vspace{3mm}

The article is based on a lecture series given by the author for the workshop ``Summer Program in PDEs" hosted by UT Austin (and conducted online) in May 2021. The author is very grateful to Philip Isett and Francesco Maggi for organizing the event. This work was supported by NSF grant DMS-1854788.

\newpage
\section{Hilbert's $19^{\text{th}}$ Problem}
Certain partial differential equations admit only $C^{\omega}$ (real analytic) solutions. An important example is the Laplace equation
\begin{equation}\label{Laplace}
\Delta u = \text{div}(\nabla u) = \delta_{ij}u_{ij} = 0
\end{equation}
from complex analysis. Another is the minimal surface equation
\begin{equation}\label{MSE}
\sqrt{1 + |\nabla u|^2}\,\text{div}\left(\frac{\nabla u}{\sqrt{1 + |\nabla u|^2}}\right) =  \left(\delta_{ij} - \frac{u_iu_j}{1 + |\nabla u|^2}\right)u_{ij} = 0,
\end{equation}
whose solutions model soap films. The solutions to these equations are $C^{\omega}$ even if their boundary values are not. Take for example the angle function $\text{Im}(\log z)$ in the half-plane $\{\text{Re}(z) > 0\} \subset \mathbb{C} \cong \mathbb{R}^2$. This function solves both equations, and is $C^{\omega}$ in $\{\text{Re}(z) > 0\}$, but is discontinuous on the boundary.

In the statement of Hilbert's $19^{\text{th}}$ problem, it is noted that equations with this remarkable property tend to arise as Euler-Lagrange equations of variational integrals of the form
\begin{equation}\label{VarInt}
 J(u) := \int F(\nabla u)\,dx,
\end{equation}
where $F$ is analytic, convex, and $\det D^2F > 0$. The Laplace equation corresponds to the choice $F(\cdot) = |\cdot|^2$, and the minimal surface equation to the choice $F(\cdot) = \sqrt{1+|\cdot|^2}$. Hilbert's $19^{\text{th}}$ problem asks whether all such Euler-Lagrange equations
\begin{equation}\label{EL}
\text{div}(\nabla F(\nabla u)) = F_{ij}(\nabla u)u_{ij} = 0
\end{equation}
admit only analytic solutions, even if the solutions have non-analytic boundary data. Henceforth we will consider this problem for functions on the unit ball $B_1 \subset \mathbb{R}^n$.

Bernstein showed in 1904 that if $n = 2$ and $u \in C^3(B_1)$ solves (\ref{EL}), then $u \in C^{\omega}(B_1)$. The regularity required on $u$ to conclude analyticity, as well as the dimension restriction, were relaxed in the following years by Lewy, Hopf, Schauder, and others (see \cite{Mor} Ch. 5.8 and the references therein). By the early 1930s, it was known that solutions to (\ref{EL}) that are in $C^{1,\,\alpha}(B_1)$ for some $\alpha > 0$ are analytic.

\begin{rem}
We say that a Lipschitz function $u$ on $B_1$ solves (\ref{EL}) in the sense of distributions if
\begin{equation}\label{WeakSolution}
\int_{B_1} \nabla F(\nabla u) \cdot \nabla \psi \,dx = 0
\end{equation}
for all $\psi \in C^1_0(B_1)$. This condition is equivalent to the statement that $u$ minimizes the integral $J$ among Lipschitz functions with the same boundary data.
\end{rem} 

The idea of the aforementioned result is as follows. First, if $\nabla u \in C^{\alpha}$, then the coefficients of the equation (\ref{EL}) are H\"{o}lder continuous, hence $u \in C^{2,\,\alpha}$ by the Schauder interior estimates (see e.g. \cite{GT} Ch. 6, or \cite{HL} Ch. 3). The coefficients of (\ref{EL}) are thus $C^{1,\,\alpha}$. We can continue applying the Schauder estimates to conclude that $u$ is smooth (more generally, that $u \in C^{k,\,\alpha}$ provided $F \in C^{k + 1},$ for any $k \geq 2$). When $F$ is analytic, the analyticity of $u$ can be shown either by carefully estimating successive derivatives of $u$ to obtain $\|D^ku\|_{L^{\infty}(B_{1/2})} \leq C^kk!$ (Bernstein's technique), or by extending to a complex domain (as Lewy, Hopf did), see also \cite{Mor} Ch. 5.8.

\begin{rem}
When $n = 1$ solutions to (\ref{EL}) are linear (hence $C^{\omega}$) regardless of the regularity of $F$. However, when $n \geq 2$, smoothness of $F$ does not guarantee analyticity of $u$ (see exercises).
\end{rem}

Although these results represented significant progress on Hilbert's problem, the existence of solutions to the Dirichlet problem
\begin{equation}\label{DP}
\begin{cases}
\text{div}(\nabla F(\nabla u)) = 0 \text{ in } B_1,\\
u|_{\partial B_1} = \varphi
\end{cases}
\end{equation}
in the class $C^{1,\,\alpha}(B_1)$ was not known. Provided for example that $\varphi \in C^2(\partial B_1)$, one can prove the existence of a Lipschitz function that solves (\ref{DP}) in the sense of distributions by minimizing the integral $J$. In the early 1930s the main problem was thus to fill the gap from Lipschitz to $C^{1,\,\alpha}$ regularity. Our first goal in these lectures will be to show how this gap was filled. 

More precisely, we will show that solutions of (\ref{EL}) satisfy estimates of the form
\begin{equation}\label{C1alpha}
\|\nabla u\|_{C^{\alpha}(B_{1/2})} \leq C\left(n,\|\nabla u\|_{L^{\infty}(B_1)},\, F\right)
\end{equation}
for some $\alpha\left(n,\|\nabla u\|_{L^{\infty}(B_1)},\, F\right) \in (0,\,1)$. To emphasize ideas, we will assume that $u$ is smooth, and establish (\ref{C1alpha}) as an a priori estimate. 

\begin{rem}
Such a priori estimates are sufficient for many purposes, for example proving the existence of classical solutions to (\ref{DP}) with sufficiently regular boundary data, when combined with appropriate tools from functional analysis.
\end{rem}

\noindent The approach to the estimate (\ref{C1alpha}) is to differentiate the equation (\ref{EL}), giving an equation in divergence form for the derivatives of $u$:
\begin{equation}\label{DEL}
\partial_i(F_{ij}(\nabla u)(\partial_ku)_j) = 0, \quad k = 1,\,...,\,n.
\end{equation}
Since we do not yet control the modulus of continuity of $\nabla u$, the idea is to treat (\ref{DEL}) as a linear, uniformly elliptic equation of the form
\begin{equation}\label{Lin}
\partial_i(a_{ij}(x)\partial_jv) = 0
\end{equation}
for $v = \partial_ku$, where the eigenvalues of $a_{ij}$ are in $[\lambda,\,\lambda^{-1}]$ for some $\lambda > 0$. Provided the estimate
\begin{equation}\label{IntEst}
\|v\|_{C^{\alpha}(B_{1/2})} \leq C(n,\,\lambda)\|v\|_{L^{\infty}(B_1)}
\end{equation}
is true for solutions to (\ref{Lin}) for some $\alpha(n,\lambda) > 0$, the key estimate (\ref{C1alpha}) follows. The estimate (\ref{IntEst}) was proven by Morrey in two dimensions in the late 1930s, and by De Giorgi \cite{DG} and Nash \cite{Na} in higher dimensions in the late 1950s. This furnished a complete solution to Hilbert's problem.

\subsection{Exercises}

\begin{enumerate}[1.]
\item Prove using the convexity of $F$ that a Lipschitz function $u$ on $B_1$ solves (\ref{EL}) in the sense of distributions if and only if it minimizes the integral $J$ (subject to its own boundary data).

\vspace{3mm}

\item Let $a_{ij}(\cdot)$ be a smoothly varying, positive definite, symmetric matrix field on $\mathbb{R}^n$. Assume that $u$ is smooth and solves
$$a_{ij}(\nabla u)u_{ij} = 0$$
in a domain $\Omega \subset \mathbb{R}^n$. Show that for any $e \in \mathbb{S}^{n-1}$, the directional derivative $u_e$ satisfies the maximum principle (that is, attains its maximum and minimum on $\partial \Omega$). In the case that $\Omega = B_1,\, \varphi \in C^2(\partial B_1)$, and $u|_{\partial B_1} = \varphi$, prove using linear functions as boundary barriers that
$$\|\nabla u\|_{L^{\infty}(B_1)} \leq C(n)\|\varphi\|_{C^2(\partial B_1)}.$$

\vspace{3mm}

\item Let $H$ be a $C^2$, even, uniformly convex function of one variable. Let $H^*$ be its Legendre transform, defined by 
$$H^*(x) = \int_0^x (H')^{-1}(s)\,ds.$$ 
Show that if $F(p,\,q) = H(p) + H(q)$, then $u(x,\,y) = H^*(x) - H^*(y)$ solves $F_{ij}(\nabla u)u_{ij} = 0.$ Using this observation, build a smooth and uniformly convex function $F$ on the plane such that the equation (\ref{EL}) has non-analytic solutions.
\end{enumerate}

\newpage
\section{Solution in Two Dimensions}
In this section we prove Morrey's estimate (\ref{IntEst}) for solutions to (\ref{Lin}) in dimension $n = 2$. We begin with a few observations that hold in any dimension. The first is that solutions to (\ref{Lin}) minimize the integral
\begin{equation}\label{LinInt}
E(v) := \int_{B_1} a_{ij}(x)v_iv_j\,dx.
\end{equation}
This follows immediately from integration by parts. Similarly, if $v$ is a subsolution to (\ref{Lin}), that is,
$$\partial_i(a_{ij}(x)v_j) \geq 0,$$
then ``downward perturbations" increase energy: $E(v - \psi) \geq E(v)$ for all non-negative functions $\psi \in C^1_0(B_1)$. It is a good exercise to show this.

 There are two important consequences. The first is the so-called Caccioppoli inequality. By choosing $v - \epsilon v \psi^2$ as a competitor for $v$, with $\psi \in C^1_0(B_1)$, and taking $\epsilon \rightarrow 0$, we arrive at
$$\int_{B_1} a_{ij}v_i(v\psi^2)_j\,dx = 0.$$
Expanding the derivative of $v\psi^2$, applying the bounds on the eigenvalues of $a_{ij}$ and using Cauchy-Schwarz gives the Caccioppoli inequality
\begin{equation}\label{Caccioppoli}
\int_{B_1} |\nabla v|^2\psi^2\,dx \leq 4\lambda^{-4}\int_{B_1} v^2|\nabla \psi|^2\,dx.
\end{equation}
This can be viewed as a ``backwards Poincar\'{e} inequality" for solutions to (\ref{Lin}). Inequality (\ref{Caccioppoli}) also holds for any nonnegative sub-solution to (\ref{Lin}), because $-v\psi^2$ is in that case a downward perturbation.

The second consequence is the maximum principle. If for some constant $h$ the set $\{v > h\}$ has a connected component that is compactly contained in $B_1$, then by replacing $v$ with $\min\{v,\,h\}$ on this component we get a competitor with lower energy $E$, a contradiction. Thus, $v$ has no interior local maxima. (The same holds for subsolutions, because the competitor is smaller.) A similar argument shows that a solution to (\ref{Lin}) has no local minima.

\begin{rem}
Another way to see that $v$ satisfies the maximum principle is to pass the derivative in the equation for $v$ to obtain
$$a_{ij}(x)v_{ij} + (\partial_ia_{ij}(x))v_j = 0,$$
and apply the maximum principle for equations in non-divergence form.
\end{rem}

Now we specialize to two dimensions. The Courant-Lebesgue lemma says that if $w$ is a function on $B_1 \subset \mathbb{R}^2$ that satisfies the maximum principle, then
\begin{equation}\label{CL}
(\text{osc}_{\partial B_r} w)^2 \leq \frac{\pi}{\log[1/(2r)]}\int_{B_{1/2}}|\nabla w|^2\,dx
\end{equation}
for $r \in (0,\,1/2)$.
Here $\text{osc}_{\Omega}w := \sup_{\Omega}w - \inf_{\Omega}w$. To prove (\ref{CL}), note first that by the fundamental theorem of calculus and Cauchy-Schwarz, we have
$$\text{osc}_{\partial B_s} w \leq \int_{\text{half of }\partial B_s} |\nabla w| \leq (\pi s)^{1/2}\left(\int_{\partial B_s} |\nabla w|^2\right)^{1/2}.$$
The inequality follows by squaring, dividing by $s$, integrating from $r$ to $1/2$, and using that $\text{osc}_{\partial B_s} w$ is non-decreasing in $s$ by the maximum principle. 

\begin{rem}
The Courant-Lebesgue lemma implies that in two dimensions, solutions to (\ref{Lin}) have a logarithmic modulus of continuity. The philosophy is that the energy $E$ is comparable to the $H^1$ norm of $v$, which in two dimensions nearly controls the modulus of continuity of $v$ by standard embeddings. The extra ingredient we use to get continuity is the maximum principle.
\end{rem}

To improve to H\"{o}lder regularity we use the scaling properties of the equation (\ref{Lin}). We show in particular that for some $\delta(\lambda) > 0$ and all $r \leq 1$ we have
\begin{equation}\label{OscDecay}
\text{osc}_{B_{\delta r}}v \leq \frac{1}{2}\text{osc}_{B_{r}}v.
\end{equation}
Iterating this inequality gives
$$\text{osc}_{B_{\delta^k}}v \leq 2^{-k}\text{osc}_{B_1}v := (\delta^k)^{\alpha}\text{osc}_{B_1}v,$$
where $\alpha$ is defined by
$$\delta^{\alpha} = 1/2.$$
This in turn implies the desired interior H\"{o}lder estimate (\ref{IntEst}) for this value of $\alpha$. We leave it as an exercise.

To prove (\ref{OscDecay}) we may assume after performing a dilation and multiplying by a constant, neither of which change the type of equation (\ref{Lin}), that $r = 1$ and that $\text{osc}_{B_1}v = \text{osc}_{\partial B_1}v = 1$. After adding a constant, which doesn't change the equation, we may assume that $0 \leq v \leq 1$. Using Courant-Lebesgue and Caccioppoli with a standard cutoff function $\psi$ that is $1$ in $B_{1/2}$, we have
$$(\text{osc}_{B_{\delta}}v)^2 \leq \frac{\pi}{\log[1/(2\delta)]}\int_{B_{1/2}} |\nabla v|^2\,dx \leq \frac{C(\lambda)}{\log[1/(2\delta)]}\int_{B_1} v^2\,dx \leq \frac{C(\lambda)}{\log[1/(2\delta)]}.$$
We arrive at the desired estimate provided $\delta(\lambda)$ is chosen small.

\begin{rem}
The Courant-Lebesgue lemma also gives a Harnack inequality for positive solutions of (\ref{Lin}) in the plane, from which the strong maximum principle and a H\"{o}lder estimate follow in a standard way. Writing $v = e^{w}$ and applying the equation for $v$, we get
$$\partial_i(a_{ij}w_j) + a_{ij}w_iw_j = 0.$$
Multiplying this by the square of a smooth standard cutoff function, integrating by parts, and using Cauchy-Schwarz, we see that
$$\int_{B_{1/2}} |\nabla w|^2 < C(\lambda).$$
Since $w$ still satisfies the maximum principle (by the monotonicity of the exponential), Courant-Lebesgue implies that
$$\text{osc}_{B_{1/4}} w \leq C(\lambda),$$
which is equivalent to the Harnack inequality
$$\sup_{B_{1/4}}v \leq C(\lambda)\inf_{B_{1/4}}v.$$
\end{rem}

\begin{rem}
In dimension $n = 2$ one can avoid using Courant-Lebesgue (and in particular the maximum principle) by using the Caccioppoli inequality more carefully. Indeed, if one takes $\psi \equiv 1$ 
in $B_{1/2}$ in the Caccioppoli inequality, and uses the invariance of (\ref{Lin}) under adding constants, one obtains
$$\int_{B_{1/2}}|\nabla v|^2\,dx \leq C(\lambda) \int_{B_1 \backslash B_{1/2}} (v-c)^2\,dx.$$
The sides of this inequality scale differently. Taking $c$ to be the average over the annulus of $v$ and applying the Poincar\'{e} inequality gives
$$\int_{B_{1/2}} |\nabla v|^2\,dx \leq C(n,\,\lambda) \int_{B_1 \backslash B_{1/2}} |\nabla v|^2\,dx,$$
that is, the mass of $|\nabla v|^2$ decays by a fixed fraction when passing from $B_1$ to $B_{1/2}$. The sides of this estimate scale the same way. By rescaling and iterating this estimate we get that
$$\int_{B_r} |\nabla v|^2\,dx \leq C(n,\,\lambda)r^{2\alpha}\int_{B_1} v^2\,dx$$
for some $\alpha(n,\,\lambda) > 0$ and all $r \leq 1/2$. In two dimensions, the previous inequality gives a $C^{\alpha}$ estimate for $v$ by Morrey space embeddings (see \cite{GT} Ch. 7 or \cite{HL} Ch. 3). An advantage of this approach is that it applies in settings in which a maximum principle is unavailable, e.g. in vector-valued problems. In contrast, in dimension $n \geq 3$, the lack of a maximum principle in vector-valued settings is fatal to regularity (see e.g. \cite{M3} and the references therein). For scalar problems in dimension $n \geq 3$, both the maximum principle and the fact that the sides of the Caccioppoli inequality scale differently play a crucial role in the proof of $C^{\alpha}$ regularity (see the next section).
\end{rem}

\begin{figure}
 \begin{center}
    \includegraphics[scale=0.7, trim={30mm 160mm 10mm 0mm}, clip]{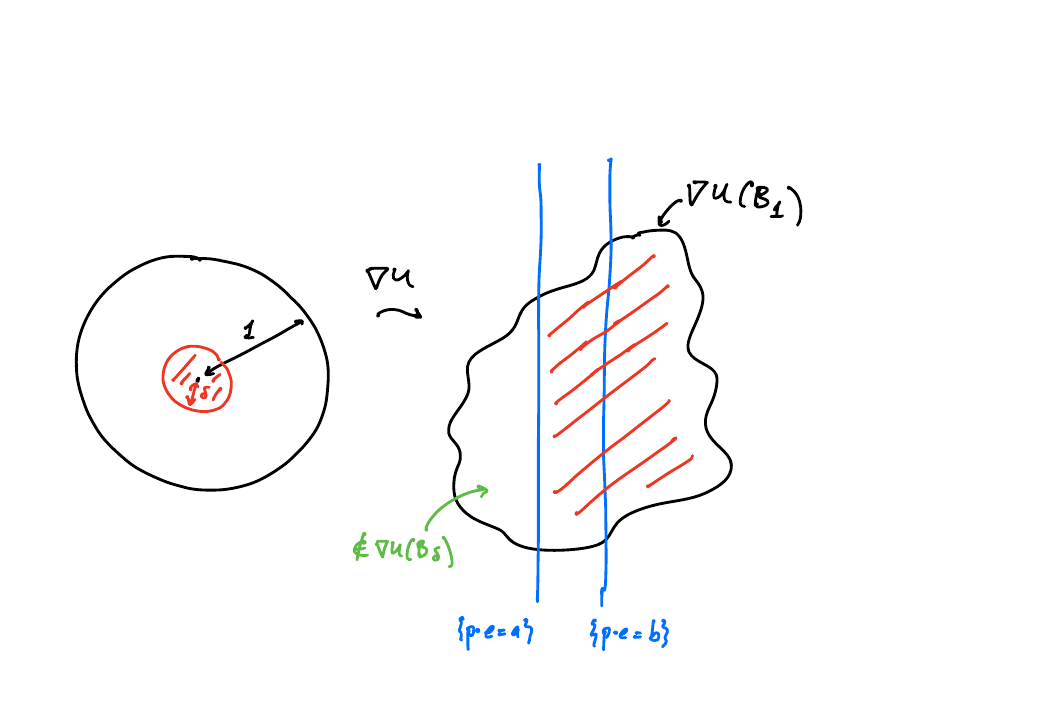}
\caption{Chopping at the gradient image with lines}
\label{LineChop}
\end{center}
\end{figure}

\begin{rem}
In the context of the original problem, where we consider Lipschitz solutions to 
$$F_{ij}(\nabla u)u_{ij} = 0,$$
the above discussion says the following. Consider the lines $\{p \cdot e = a\}$ and $\{p \cdot e = b\}$, with $a < b$, in $\mathbb{R}^2$. Then as $r \rightarrow 0^+$, the sets
$\nabla u(B_r)$ must localize to one of the half-spaces $\{p \cdot e \leq b\}$ or $\{p \cdot e \geq a\}$, that is, we can ``chop at the gradient image of $u$" with level sets of linear functions (see Figure \ref{LineChop}). Indeed, if $\nabla u(B_{\delta})$ crosses the strip $\{p \cdot e \in [a,\,b]\}$ for $\delta$ small, then by the maximum principle it does so on all circles $\partial B_r$ with $r \geq \delta$. The Courant-Lebesgue lemma says that the $L^2$ norm of $D^2u$ is huge in $B_{1/2}$, but this violates the Caccioppoli inequality. By varying the choice of direction $e$ and repeating this argument, we see that $\nabla u(B_r)$ localizes to a point as $r \rightarrow 0^+$, i.e. $C^1$ regularity.
\end{rem}

\begin{rem}
The approach we outlined to solving Hilbert's problem involved differentiating the Euler-Lagrange equation (\ref{EL}). In two dimensions it turns out that this is not necessary. Assume that $w$ solves a linear uniformly elliptic equation of the form
$$a_{ij}(x)w_{ij} = 0$$
in $\mathbb{R}^2$, where the eigenvalues of the coefficient matrix are in $[\lambda,\,\lambda^{-1}]$ for some $\lambda > 0$. Then the eigenvalues of $D^2w$ have opposite sign and comparable absolute value, i.e.
$$|D^2w|^2 \leq -C(\lambda)\det D^2w.$$
The map $\nabla w$ is thus a quasi-conformal map, a generalization of a holomorphic map that infinitesimally takes disks to ellipses with bounded eccentricity. Such maps are well-studied, and satisfy interior $C^{\alpha}$ estimates, implying interior $C^{1,\,\alpha}$ regularity for solutions to such equations in two dimensions (see \cite{GT} Ch. 12). In higher dimensions, this result is false. In \cite{Saf} Safonov constructed, for each $\alpha \in (0,\,1)$, functions on $\mathbb{R}^3$ that are homogeneous of degree $\alpha$ and solve (in the viscosity sense) linear uniformly elliptic equations in non-divergence form. It is not a coincidence that such examples were not constructed with $\alpha = 1$, and we will revisit this in a later section.
\end{rem}

\subsection{Exercises}
\begin{enumerate}[1.]
\item Construct an example of a function $v$ on $B_1 \subset \mathbb{R}^3$ such that $\text{osc}_{\partial B_r} v$ is non-decreasing in $r$ and $\int_{B_1} |\nabla v|^2\,dx < \infty$, but $v$ is discontinuous at the origin. Hint: Choose any zero-homogeneous function that is smooth and non-constant on $\mathbb{S}^2$.

\vspace{3mm}

\item Show that the functions $u(x) = |x|^{\alpha}g(x/|x|)$, where $g$ is a non-constant eigenfunction of $\Delta_{\mathbb{S}^{n-1}}$ and $\alpha > 0$, solve linear uniformly elliptic equations in divergence form. (Hint: take $a_{ij} = \delta_{ij} + \mu\frac{x_ix_j}{|x|^2}$ for appropriate $\mu > -1$.)

\vspace{3mm}

\item Assume that $a_{ij}(x)u_{ij} = 0$ in $B_1 \subset \mathbb{R}^2$ with $a_{ij}$ uniformly elliptic. By differentiating the equation once, and using the original equation, show that $u_1$ satisfies the maximum principle. (Hint: $u_{22}$ can be written as a linear combination of derivatives of $u_1$.)

Show next that $\partial_1u$ solves a linear uniformly elliptic equation in divergence form. Conclude from the discussion in this section that $u$ enjoys $C^{1,\,\alpha}$ estimates independent of the regularity of the coefficients $a_{ij}$. (Hint: up to dividing by $a_{22}$, the equation can be written $a_{11}u_{11} + 2a_{12}u_{12} + u_{22} = 0$.)
\end{enumerate}

\newpage
\section{Solution in Higher Dimensions}
In this section we outline De Giorgi's proof of the estimate (\ref{IntEst}) in higher dimensions. His approach was inspired by the regularity theory for minimal surfaces, in particular the ``density estimate," which says that each side of a minimal hypersurface fills a nontrivial fraction of any (extrinsic) ball centered on the surface (there are no ``spikes"). The analogous result in the function case is the so-called $L^2-L^{\infty}$ estimate:

\begin{thm}\label{NoSpikes}
Assume that $\partial_i(a_{ij}(x)\partial_j v) \geq 0$ in $B_2 \subset \mathbb{R}^n$. Then
\begin{equation}\label{L2LInfinity}
\sup_{B_1} v \leq C(n,\,\lambda)\|v_+\|_{L^2(B_2)}.
\end{equation}
\end{thm}

\noindent That is, sub-solutions have no interior upward spikes. Here $v_+ := \max\{v,\,0\}$.

We sketch the proof. It is a good exercise to show that for any increasing convex function $G$ of one variable, we have $\partial_i(a_{ij}(x)\partial_j(G(v))) \geq 0$. In particular, for any $\kappa \in \mathbb{R}$, the function $v_{\kappa} := (v-\kappa)_+$ satisfies the Caccioppoli inequality. Using Cauchy-Schwarz and Caccioppoli, we have for any cutoff function $\psi \in C^{\infty}_0(B_2)$ that
$$\int |\nabla (v_{\kappa}^2\psi^2)| \leq C(\lambda)\left(\int v_{\kappa}^2\psi^2\right)^{1/2}\left(\int v_{\kappa}^2|\nabla \psi|^2\right)^{1/2}.$$
Applying the Sobolev inequality to the left side we arrive at
$$\left(\int (v_{\kappa}^2\psi^2)^{\frac{n}{n-1}}\right)^{\frac{n-1}{n}} \leq C(n,\,\lambda)\left(\int v_{\kappa}^2\psi^2\right)^{1/2}\left(\int v_{\kappa}^2|\nabla \psi|^2\right)^{1/2}.$$
Combining H\"{o}lder's inequality and the previous estimate we get
$$\int v_{\kappa}^2\psi^2 \leq C(n,\,\lambda)\left(\int v_{\kappa}^2|\nabla \psi|^2\right)|\{v_{\kappa}^2\psi^2 > 0\}|^{\frac{2}{n}}.$$
Let $0 \leq \tau < \kappa \leq 1$. Assume that $\psi$ is a standard cutoff that is $1$ in $B_{2-\kappa}$ and $0$ outside of $B_{2-\tau}$. Then the previous inequality implies that
$$\int_{B_{2-\kappa}} v_{\kappa}^2 \leq \frac{C}{(\kappa-\tau)^2}\left(\int_{B_{2-\tau}} v_{\tau}^2\right)|\{v_{\tau} \geq \kappa - \tau\} \cap B_{2-\tau}|^{2/n}.$$
Denoting 
$$\int_{B_{2-s}} v_s^2 := V(s)$$ 
and applying the Chebyshev inequality to the last term on the right side of the previous inequality, we obtain
\begin{equation}\label{ScalingClass}
V(\kappa) \leq \frac{C}{(\kappa - \tau)^{2 + 4/n}}V^{1 + \frac{2}{n}}(\tau).
\end{equation}
Inequality (\ref{ScalingClass}) implies that $V(1) = 0$ provided $V(0)$ is sufficiently small (see exercises). (The key point is the difference in powers of $V$ appearing in (\ref{ScalingClass}). This arises from the use of the Sobolev inequality, which competes, in a scaling-invariant way, with the non-scaling-invariant Caccioppoli inequality.) In particular, $v \leq 1$ in $B_1$ provided $\|v_+\|_{L^2(B_2)}$ is sufficiently small depending on $n,\,\lambda$. Using the invariance of the equation (\ref{Lin}) under multiplication by constants, Theorem \ref{NoSpikes} follows.

Theorem \ref{NoSpikes} can be used to prove the following useful oscillation decay result:

\begin{prop}\label{nDOsc}
Assume that $\partial_i(a_{ij}(x)\partial_jv) \geq 0$ in $B_1 \subset \mathbb{R}^n$. For $\delta > 0$, there exists $\epsilon(n,\,\lambda,\,\delta) > 0$ such that if $|\{v_+ = 0\} \cap B_1| \geq \delta |B_1|,$ then
$$\sup_{B_{1/2}}v \leq (1-\epsilon)\sup_{B_1}v_+.$$
\end{prop}

\begin{figure}
 \begin{center}
    \includegraphics[scale=0.6, trim={30mm 160mm 10mm 0mm}, clip]{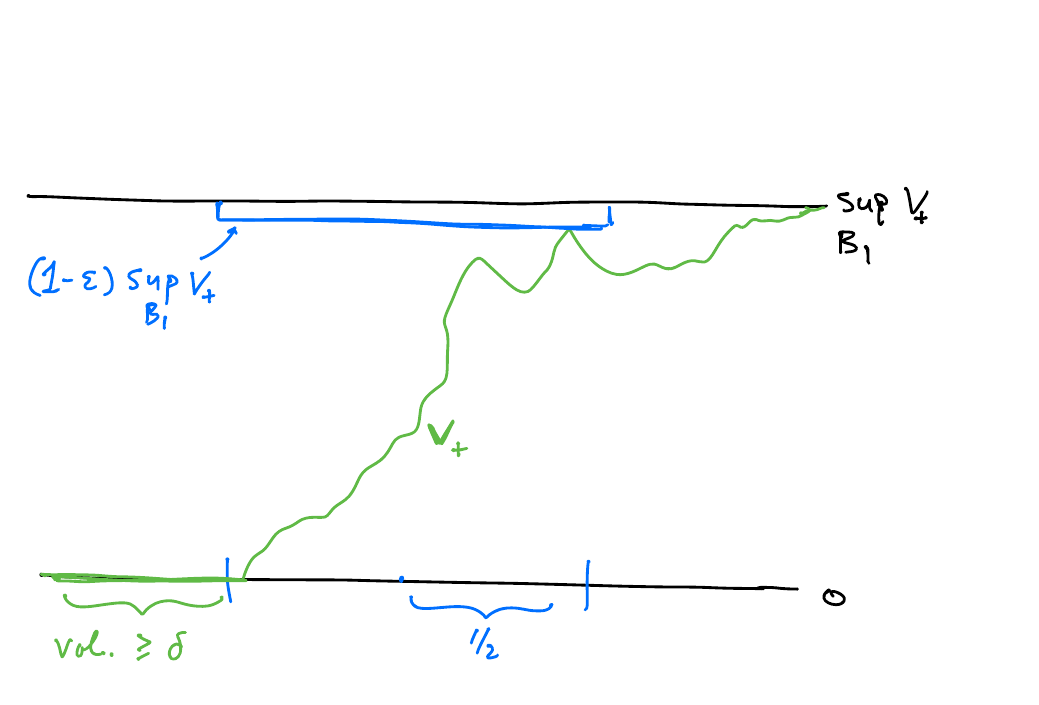}
\caption{Strict interior separation from maximum}
\label{OscDrop}
\end{center}
\end{figure}

\noindent Proposition \ref{nDOsc} says that sub-solutions that are not close to their maximum nearly everywhere must separate from their maximum when we step away from the boundary (see Figure \ref{OscDrop}). This is a quantitative version of the maximum principle, which says that sub-solutions don't have interior maxima.

We will prove a slightly weaker version of Proposition \ref{nDOsc} to minimize technicalities and emphasize ideas. Namely, we will assume that $v$ is a sub-solution in $B_2$ and that $|\{v_+ = 0\} \cap B_1| \geq \delta |B_1|$, and we will show that $\sup_{B_{1/2}} v \leq (1-\epsilon) \sup_{B_2}v^+$. To that end assume that $\sup_{B_2}v_+ = 1$ and let
$$W(s) := \frac{|\{v \leq s\} \cap B_1|}{|B_1|}.$$
Then $W \in [0,\,1]$ is nondecreasing, and we will show that it increases to $1$ quantitatively. Let $0 \leq s < t \leq 1$, and let
$$w = \frac{(v-s)_+}{1-s}$$
so that $0 \leq w \leq 1$. Let $\bar{w} = \min\{w,\, (t-s)/(1-s)\}$. By Cauchy-Schwarz and Caccioppoli we have
$$\int_{B_1} |\nabla \bar{w}| \leq \left(\int_{B_1} |\nabla w|^2\right)^{1/2}|B_1|^{1/2}(W(t) - W(s))^{1/2} \leq C(n,\,\lambda)(W(t)-W(s))^{1/2}.$$
On the other hand, we have by the Poincar\'{e} inequality that
$$\int_{B_1} |\nabla \bar{w}| \geq c(n)\int_{B_1} |\bar{w} - \text{avg.}_{B_1}\bar{w}| \geq c(n)W(s)(1-W(t))\frac{t-s}{1-s}.$$
Putting the previous two inequalities together we have
$$W(s)\left(1 + c(n,\,\lambda)\frac{(t-s)^2}{(1-s)^2}W(s)(1-W(t))^2\right) \leq W(t).$$
Geometrically, this inequality says that $H^1$ functions ``pay in measure" to pass from one height to another.

Now assume that $W(0) \geq \delta > 0$ and let $c$ denote a small constant depending only on $n,\,\lambda,\,\delta$ which may change from line to line. Let $b_k = W(1-2^{-k})$. The previous inequality gives
$$b_k(1 + c(1-b_{k+1})^2) \leq b_{k+1}.$$
Letting $a_k = 1-b_k$ this becomes
$$a_{k+1} + ca_{k+1}^2 \leq a_k.$$
Since $a_k \leq 1$ this implies that
$$a_{k+1} \leq a_k - ca_k^2.$$
It follows that $a_k \leq \frac{1}{1+ ck}$ (see exercises). In particular, there is some $\epsilon(n,\,\lambda,\,\delta) > 0$ such that the sub-solution 
$\bar{v} = \frac{[v-(1-2\epsilon)]_+}{2\epsilon}$ has $L^2$ norm small enough in $B_1$ that Theorem \ref{NoSpikes} gives $|\bar{v}| \leq 1/2$ in $B_{1/2}$, i.e. $v \leq 1-\epsilon$ in $B_{1/2}$.

\subsection{Exercises}
\begin{enumerate}[1.]
\item Prove that if $a_{k+1} \leq C^ka_k^{1+\gamma}$ for some $C,\,\gamma > 0$, then provided $a_0$ is sufficiently small depending on $C$ and $\gamma$, we have
$a_k \rightarrow 0$. Show using (\ref{ScalingClass}) that $a_k := V(1-2^{-k})$ satisfies such a relation.

\vspace{3mm}

\item Prove that if $0 \leq a_{k+1} \leq a_k - ca_k^2$ for some $c > 0$ small, and $a_0 \in [0,\,1]$, then $a_k \leq \frac{1}{1+ck}$ for all $k \geq 0$.

\vspace{3mm}

\item Prove the estimate (\ref{IntEst}) for solutions to (\ref{Lin}) in all dimensions using Proposition \ref{nDOsc}. (Hint: It suffices to prove the oscillation decay estimate
$$\text{osc}_{B_{1/2}}v \leq (1-\epsilon(n,\,\lambda))\text{osc}_{B_1}v.$$
One may assume after adding a constant and multiplying by a constant that $\inf_{B_1} v = -1$ and $\sup_{B_1}v = 1$. One of $|\{v \leq 0\}|$ or $|\{v \geq 0\}|$ is at least $|B_1|/2$. Use Proposition \ref{nDOsc} to conclude.)

\end{enumerate}

\newpage
\section{Degenerate Convex Functionals in Two Dimensions}
In this section we begin to discuss the regularity of Lipschitz minimizers of (\ref{VarInt}) in the case that $F$ is not smooth and uniformly convex. One important example is the $p$-Laplace energy density
$$F(\cdot) = |\cdot|^p,$$
with $p > 1$ and $p \neq 2$. Another important example, which arises in models of traffic congestion (see e.g. \cite{CF}), is
$$F(\cdot) = (|\cdot|-1)_+^2.$$ 
The latter integrand vanishes in $B_1$, so any $1$-Lipschitz function is a minimizer of the corresponding functional. Such ``degenerate convex" Lagrangians also arise in models of crystal surfaces (see for example \cite{KOS}, \cite{DMMN}).

The existence of Lipschitz minimizers of (\ref{VarInt}) with sufficiently regular boundary data is not hard to show. A natural question is which conditions on $F$ guarantee that Lipschitz minimizers (equivalently, Lipschitz solutions to (\ref{EL})) are $C^1$. Obtaining such a result would be useful for understanding finer properties of solutions. For example, if $u \in C^1$ and $\nabla u(x_0)$ lies in a region where $F$ is smooth and uniformly convex, then the classical theory would imply that $u$ is smooth nearby $x_0$.

\begin{rem}
The {\it local} Lipschitz regularity of minimizers of (\ref{VarInt}) for Lagrangians that satisfy various growth conditions at infinity is a delicate and active research topic, with important contributions by many authors- see the survey \cite{Min} and the references therein. 
\end{rem}

If the graph of $F$ contains line segments, it is straightforward to construct Lipschitz minimizers that are not $C^1$. Indeed, after subtracting a linear function from $F$ (which does not change the equation) we may assume that the minimum set of $F$ contains a line segment from $ae$ to $be$ for some $a < b$ and $e \in \mathbb{S}^{n-1}$. Then any function $u(x) = g(x \cdot e)$ with $g' \in [a,\,b]$ is a minimizer.

Another observation is that the Legendre transform $F^*$ of $F$ solves 
$$\text{div}(\nabla F(\nabla F^*)) = F_{ij}(\nabla F^*)F^*_{ij} = n,$$
which resembles the Euler-Lagrange equation (\ref{EL}) but with nonzero constant right-hand side. The function $F^*$ is $C^1$ if and only if $F$ is strictly convex. 

These observations motivate the question of whether minimizers have the same regularity as $F^*$, and in particular, whether minimizers are $C^1$ when $F$ is strictly convex. The answer turns out to be ``no" in general (shown recently in \cite{M2}), but ``yes" in special cases. 

Here and below we assume that $F$ is convex on $\mathbb{R}^n$, and that off of some compact degeneracy set $K$ the function $F$ is smooth with $\det D^2F > 0$. In this section we will discuss the result of De Silva-Savin \cite{DS} that Lipschitz minimizers of (\ref{VarInt}) are $C^1$ when $n = 2$ and $K$ is finite. 

\begin{rem}
The $C^1$ regularity of Lipschitz minimizers in any dimension in the case $K = \emptyset$ is the De Giorgi-Nash theorem. Lipschitz minimizers are also $C^1$ in any dimension when $K$ is a single point; this case can be treated using ideas from the theory of the $p$-Laplace equation, which was studied by many authors including Ural'tseva \cite{Ura}, Uhlenbeck \cite{Uhl}, Evans \cite{E}, Lewis \cite{L}, Tolksdorff \cite{T}, and others. We will discuss a more general result in the next section. 
\end{rem}

Our discussion in this section and the next section will be qualitative rather than quantitative, and we will also freely differentiate the equation, in order to emphasize ideas.

We first examine what can be done with the tools developed for the nondegenerate case $K = \emptyset$. We used that linear functions of $\nabla u$ solve elliptic equations. An important observation is that arbitrary convex functions of $\nabla u$ are sub-solutions of elliptic equations, and thus take their maxima on the boundary of their domains of definition. Indeed, we compute
\begin{equation}\label{EtaEqn}
\partial_i(F_{ij}(\nabla u)\partial_j(\eta(\nabla u))) = F_{ij}u_{jk}\eta_{kl}u_{li} = \sum_{k = 1}^n \lambda_k D^2F|_{\nabla u}(\nabla u_k,\,\nabla u_k)
\end{equation}
in coordinates where $\eta_{kl} = \lambda_k\delta_{kl}$. If $\eta$ is convex then $\lambda_k \geq 0$, so this expression is nonnegative by the convexity of $F$. Furthermore, if $\eta \geq 0$ and $\eta$ vanishes on a region containing $K$, then $\eta(\nabla u)$ is a sub-solution to a {\it uniformly} elliptic equation, because the coefficients $F_{ij}(\nabla u)$ play no role in the equation for $\eta(\nabla u)$ on the set $\{\eta(\nabla u) = 0\}$. In particular, $\eta(\nabla u)$ satisfies the Caccioppoli inequality.

Take for example the function $\eta(p) = (p \cdot e - a)_+$, with $a$ and $p$ chosen such that 
$$K \subset \{p \cdot e < a\}.$$ 
The maximum principle applied to linear functions of $\nabla u$ says that if $\eta(\nabla u) = 0$ at some point on $\partial B_{\delta}$, then the same holds on all larger spheres. Similarly, if $h > 0$ and $\eta(\nabla u) \geq h$ at some point on $\partial B_{\delta}$ (that is, $\nabla u(B_{\delta})$ intersects $\{p \cdot e \geq a + h\}$), the same holds on all larger spheres. We now restrict our attention to the case $n = 2$. For $\delta$ small, the Courant-Lebesgue lemma says that the $H^1$ norm of $\eta(\nabla u)$ is very large. However, the Caccioppoli inequality for $\eta(\nabla u)$ implies that this norm is bounded depending on $\|\nabla u\|_{L^{\infty}(B_1)}$ and the properties of $F$ away from $K$. For $\delta$ sufficiently small these inequalities are incompatible, hence $\nabla u(B_{\delta})$ is contained in one of the half-spaces $\{p \cdot e > a\}$ or $\{p \cdot e < a + h\}$. In the former case, the equation for $u$ is uniformly elliptic in $B_{\delta}$ and the classical theory can be applied. In the latter case, we can repeat the argument for the function $\delta^{-1}u(\delta \cdot)$, which also solves (\ref{EL}). Iterating this procedure, we see that as $r \rightarrow 0$, the set $\nabla u(B_r)$ localizes either to a point outside the convex hull of $K$, or to the convex hull of $K$, in dimension $n = 2$. (The same is true in higher dimensions, but the proof is more involved; see the next section).

The issue is thus to decide whether $\nabla u(B_r)$ localizes beyond the convex hull of $K$ as $r \rightarrow 0$. Work of De Silva-Savin \cite{DS} gives a way of localizing $\nabla u$ to a connected component of $K$ in dimension $n = 2$. The key is that certain non-convex functions of $\nabla u$ are sub-solutions to uniformly elliptic equations. Indeed, consider the calculation (\ref{EtaEqn}) above. Assume we have chosen a point $x$ such that $\nabla u(x) \notin K$, so $D^2F(\nabla u(x))$ has eigenvalues bounded between positive constants that we control. Below $c,\,C$ will denote small and large positive numbers depending on these constants. Up to permuting coordinates may assume that $\lambda_1 \leq \lambda_2 \leq ... \leq \lambda_n$. If $\lambda_2 > 0$ and $\lambda_1 > -\epsilon \lambda_2$, then we have
$$\sum_{k = 1}^n \lambda_k D^2F|_{\nabla u}(\nabla u_k,\,\nabla u_k) \geq c\lambda_2\left(\sum_{(i,\,j) \neq (1,\,1)} u_{ij}^2 - C\epsilon \sum_{i \geq 2} u_{1i}^2 - C\epsilon u_{11}^2\right).$$
For $\epsilon$ small the second term on the right side can be absorbed by the first. Using the equation $F_{11}(\nabla u)u_{11} = -\sum_{(i,\,j) \neq (1,\,1)} F_{ij}(\nabla u)u_{ij}$ and taking $\epsilon$ small we see that the last term can also be absorbed by the first, hence the expression is nonnegative. Thus, functions of $\nabla u$ with all positive Hessian eigenvalues except for one slightly negative one are sub-solutions at points where $\nabla u$ lies in the non-degenerate region. We can thus hope to localize $\nabla u(B_r)$ using nonnegative functions of $\nabla u$ that vanish on $K$, and have level sets that bend ``away from $K$" in exactly one direction and ``towards $K$" in the remaining directions. 

\begin{figure}
 \begin{center}
    \includegraphics[scale=0.6, trim={30mm 160mm 10mm 0mm}, clip]{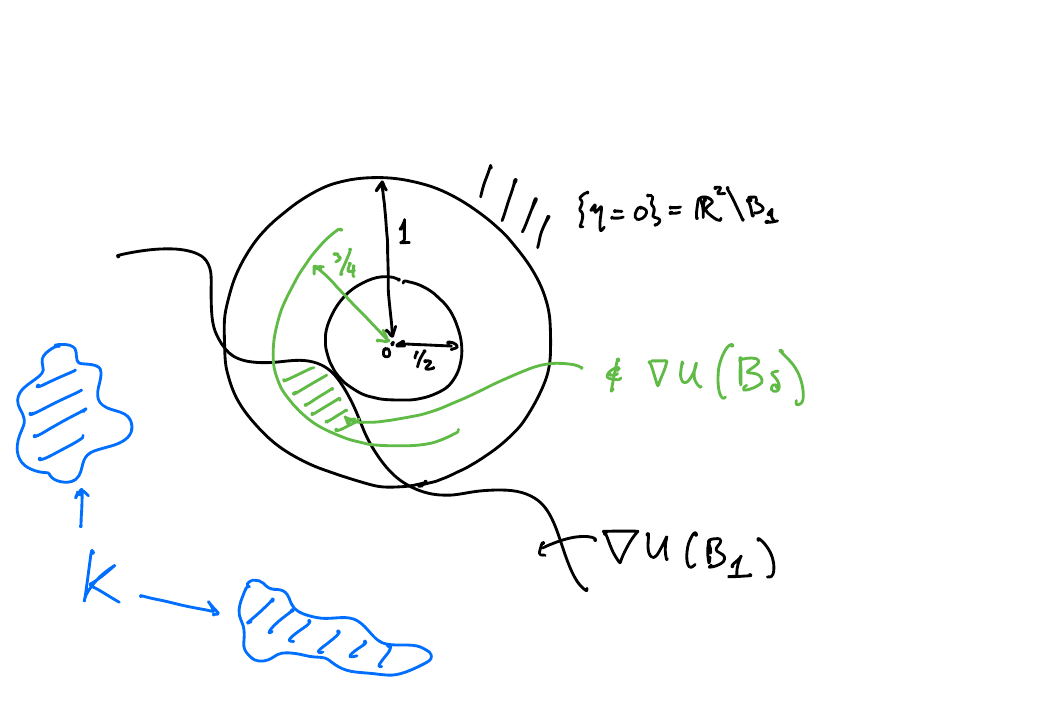}
\caption{Chopping with circles}
\label{CircleChop}
\end{center}
\end{figure}

To see this in action in two dimensions, assume for example that $\partial B_{1/2}$ touches $\nabla u(B_1)$ from the exterior, and that $K$ lies outside $\overline{B_1}$. Thus, at points where $\nabla u$ lies in $B_1$, the equation (\ref{EL}) is nondegenerate. Let $M > 0$ and let 
$$\eta(p) = [|p|^{-M} - 1]_+.$$ 
The positive (radial) Hessian eigenvalue of $\eta$ in $B_1$ is $M(M+1)|p|^{-M-2}$, while the negative Hessian eigenvalue (tangential to circles) is $-M|p|^{-M-2}$. Thus, for $M$ sufficiently large, $\eta(\nabla u)$ is a nonnegative sub-solution to a uniformly elliptic equation. If $\eta(\nabla u) = 0$ somewhere on $\partial B_{\delta}$, this remains the case on larger circles by the maximum principle applied to the convex function $|\cdot|$ of $\nabla u$. If in addition $\eta(\nabla u) > \eta(3e_1/4)$ somewhere on $\partial B_{\delta}$, this remains the case on larger circles by the fact that $\eta(\nabla u)$ is a sub-solution to an elliptic equation. Provided $\delta$ is small, Courant-Lebesgue implies that the $H^1$ norm of $\eta(\nabla u)$ is large, but as above this violates the Caccioppoli inequality. We conclude that $\nabla u(B_{\delta})$ is either contained in $B_1$, in which case the equation is non-degenerate in $B_{\delta}$, or $\nabla u(B_{\delta})$ is outside of $B_{3/4}$, i.e. we have chopped at the gradient image with a circle (Figure \ref{CircleChop}). By chopping with circles of various sizes and locations, we see that $\nabla u(B_r)$ localizes as $r \rightarrow 0^+$ either to a connected component of $K$, or to a point outside of $K$. In particular, if $K$ is finite and $n = 2$, then $u \in C^1$.

\subsection{Exercises}
\begin{enumerate}[1.]
\item Assume that $\Delta u = 0$. Let $\eta$ be a function on $\mathbb{R}^n$ such that $D^2\eta$ has eigenvalues $\lambda_1 \leq ... \leq \lambda_n$. Prove that  if $\lambda_{2} > 0$ and $\lambda_1 \geq -\frac{1}{n-1}\lambda_2$, then $\eta(\nabla u)$ is subharmonic.
\end{enumerate}

\newpage
\section{Degenerate Convex Functionals in Higher Dimensions}
In this section we consider degenerate convex functionals in dimension $n \geq 3$. Assume as in the previous section that $u$ is a Lipschitz minimizer of (\ref{VarInt}), where $F$ is convex, and away from a compact degeneracy set $K$ the function $F$ is smooth and satisfies $\det D^2F > 0$.

We first discuss the fact that in any dimension, the sets $\nabla u(B_r)$ localize as $r \rightarrow 0^+$ to either a point outside the convex hull of $K$, or to the convex hull of $K$ (see \cite{CF}, \cite{M2}). In dimension $n \geq 3$, proving this is more subtle than in two dimensions, where the Courant-Lebesgue lemma could be used. The starting point is the same: assume that $K$ is contained in a half-space $\{p \cdot e < a\}$, and consider the sub-solution $\eta(\nabla u) = (u_e - a)_+$. If the measure of the set of points $x$ such that $\eta(\nabla u(x)) = 0$ is bounded away from zero, then we can apply Proposition \ref{nDOsc} to conclude that $\nabla u(B_{1/2})$ is contained in a region that is quantitatively smaller than $\nabla u(B_1)$, namely
$$\nabla u(B_{1/2}) \subset \{p \cdot e \leq a + (1-\epsilon)(\sup_{B_1}\eta(\nabla u))\}.$$ 
The alternative is that $\nabla u$ lies in the half-space $\{p \cdot e \geq a\}$ (away from $K$) in the vast majority of $B_1$, which morally means that the equation is already uniformly elliptic. However, this needs to be made precise. 

To that end we invoke a result of Savin \cite{Sav}, which says that if the equation (\ref{EL}) holds in $B_1$ and $u$ is sufficiently close in $L^{\infty}(B_1)$ to a linear function $L$ such that $\nabla L$ is outside of $K$, then $u$ is smooth in $B_{1/2}$ and $\nabla u$ is very close in $L^{\infty}(B_{1/2})$ to $\nabla L$. Now, the argument goes as follows: if $\nabla u(B_1)$ contains points outside the convex hull of $K$, we can choose a direction $e$ and a value $a$ such that $\{p \cdot e \geq a\} \cap \nabla u(B_1)$ has tiny diameter and lies away from $K$. The above dichotomy argument says that either $\nabla u(B_{1/2})$ is contained in $\nabla u(B_1)$ with a piece removed, or $\nabla u$ is extremely close (on average) to a point that lies outside of $K$. In the latter case we can say that $u$ is very close in $L^{\infty}$ to a linear function whose gradient lies in the non-degeneracy region for $F$, hence Savin's result applies and we are done. If the former case happens, we repeat the argument for the rescaling $2u(\cdot/2)$, which solves (\ref{EL}). Iterating the argument, we have that either $\nabla u(B_r)$ tends to a point outside the convex hull of $K$ as $r \rightarrow 0$, or if the first case in the dichotomy continues happening we have that $\nabla u(B_r)$ gets as close as we like to the convex hull of $K$.

\begin{figure}
 \begin{center}
    \includegraphics[scale=0.9, trim={10mm 180mm 10mm 0mm}, clip]{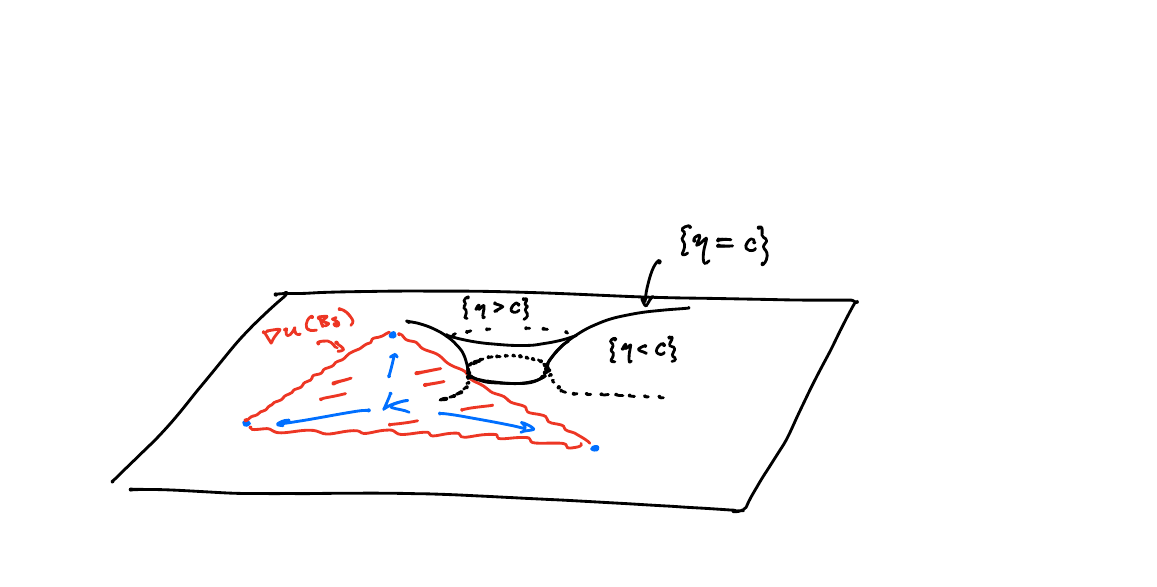}
\caption{The case that $K$ consists of $3$ points}
\label{3PointChop}
\end{center}
\end{figure}

Thus, the issue in any dimension is to localize the gradient of $u$ beyond the convex hull of $K$. In the previous section we outlined a strategy based on building non-convex functions $\eta$ of $\nabla u$ that are sub-solutions to the linearized equation. More precisely, we seek functions $\eta$ that are nonnegative and vanish on $K$, which have $n-1$ positive Hessian eigenvalues, and possibly one small negative Hessian eigenvalue. The level sets of such functions can only bend ``away from $K$'' in {\it one} direction, so one cannot for example chop from the outside using spheres in $\mathbb{R}^3$ (in contrast with the two-dimensional case, where chopping from the outside of $K$ with circles could be done). On the other hand, if $K$ has two-dimensional convex hull, the previous discussion reduces the problem to considering solutions whose gradients are very close to a two-dimensional subspace of $\mathbb{R}^n$. One can then chop with hypersurfaces that bend away from $K$ in only one direction to localize to connected components of $K$, as in the two-dimensional case (see Figure \ref{3PointChop}). Heuristically, localization to the convex hull of $K$ reduces the problem to the two-dimensional case. In particular, if $K$ is finite and contained in a $2$-plane, e.g. consists of three or fewer points, then $u \in C^1$ regardless of dimension. This was proven in \cite{M2}.

It is natural to ask for $C^1$ regularity results with less restrictive hypotheses. In view of the result in dimension $n = 2$, a natural guess is that $K$ having codimension two, or at least $K$ being finite, would suffice. However, such results are false without imposing additional structure on $F$. Namely, there are interesting counterexamples and conjectured counterexamples:

\begin{enumerate}[i.)]
\item $F$ being strictly convex doesn't imply $u \in C^1$, at least in dimension $n \geq 4$,
\item $K$ having codimension two doesn't imply $u \in C^1$, at least in dimension $n \geq 4$,
\item It seems likely that $K$ being finite doesn't imply that $u \in C^1$, in dimensions $n \geq 3$.
\end{enumerate}

\noindent In \cite{M2} a Lipschitz but non-$C^1$ minimizer to a functional of the form $\int F(\nabla u)$ is constructed in dimension $n = 4$, where $F$ is uniformly convex and $K = \mathbb{S}^1 \times \mathbb{S}^1$, which proves assertions i) and ii). We'll get to the third point below.

To build a counterexample, the idea is to start with a one-homogeneous function $u$ that is saddle-shaped: $D^2u$ is indefinite. Such a function is invariant under the rescalings that preserve the equation (\ref{EL}), and having indefinite Hessian means that $u$ solves some elliptic equation. 
An extremely useful way to build such a function is through a correspondence between one-homogeneous functions and certain singular hypersurfaces in $\mathbb{R}^n$. A one-homogeneous function $u$ gives rise to a hypersurface $\nabla u(\mathbb{S}^{n-1}) = \nabla u(\mathbb{R}^n \backslash \{0\})$, known as the ``hedgehog" of $u$. The unit normal $\nu$ to the hedgehog of $u$ and the gradient of $u$ are related by
\begin{equation}\label{HedgehogCorrespondence}
\nu(\nabla u(x)) = x
\end{equation}
for $x \in \mathbb{S}^{n-1}$, at least where $\nabla u(\mathbb{S}^{n-1})$ is smooth (see the exercises). Conversely, any hypersurface with injective Gauss map (understood in a certain generalized sense) is the hedgehog of its support function. See e.g. \cite{MM2} for a deeper discussion of hedgehog theory.

Differentiating the relation (\ref{HedgehogCorrespondence}) gives
$$II(\nabla u(x)) = (D^2u)^{-1}(x),$$
where $II$ denotes the second fundamental form of the hedgehog. Thus, choosing a candidate for a singular minimizer is equivalent to finding a hypersurface with injective Gauss map that is saddle-shaped (a ``hyperbolic hedgehog"), and taking its support function. 

Once $u$ is chosen, the game is to build the Lagrangian $F$, which by (\ref{EL}) must solve 
$$F_{ij}(\nabla u) u_{ij} = II^{ij}(\nabla u)F_{ij}(\nabla u) = 0.$$
This can be viewed as a linear equation of hyperbolic type for $F$ on the hedgehog. The challenge is to build a function $F$ that is globally convex, solves the above PDE on the hedgehog, and has a small degeneracy set $K$, most likely where the hedgehog is singular.

The example from \cite{M2} in four dimensions is
\begin{equation}\label{4DEx}
u(z_1,\,z_2) = \frac{1}{\sqrt{2}}\frac{|z_1|^2-|z_2|^2}{\sqrt{|z_1|^2 + |z_2|^2}},
\end{equation}
with $z_i \in \mathbb{C}$.
The gradient image of $u$ consists of two smooth saddle-shaped components, one in $\{|z_1| > |z_2|\}$ and the other its reflection over $\{|z_1| = |z_2|\}$, that meet on the Clifford torus $\mathbb{S}^1 \times \mathbb{S}^1$ where $\nabla u(\mathbb{S}^{3})$ is singular. The hedgehog of $u$ has enough symmetry that one can build the integrand $F$ using ODE and extension techniques. It turns out that $D^2F \geq cI$ for some $c > 0$, and that $F$ is smooth away from the Clifford torus $K := \mathbb{S}^1 \times \mathbb{S}^1$, with exactly one eigenvalue of $D^2F$ tending to infinity on $K$.

How does this approach fare in three dimensions? A first obstruction is that there are no one-homogeneous solutions to uniformly elliptic equations of the form $a_{ij}(x)u_{ij} = 0$ in three dimensions (see the paper \cite{HNY} of Han-Nadirashvili-Yuan and the exercises below). In contrast, the example (\ref{4DEx}) above does solve a linear uniformly elliptic equation in nondivergence form. The obstruction in three dimensions is, roughly speaking, that one-homogeneous functions of three variables are really two-dimensional (they have only two nonzero Hessian eigenvalues). In particular, if such a function solves an elliptic PDE, then its gradient satisfies the maximum principle (see the exercises in Section 3). Since the gradient is zero-homogeneous it is a function on the sphere (a compact manifold) and it is thus constant. The exercises outline a rigorous proof of this result, following \cite{HNY}.

\begin{figure}
 \begin{center}
    \includegraphics[scale=0.4, trim={10mm 70mm 10mm 0mm}, clip]{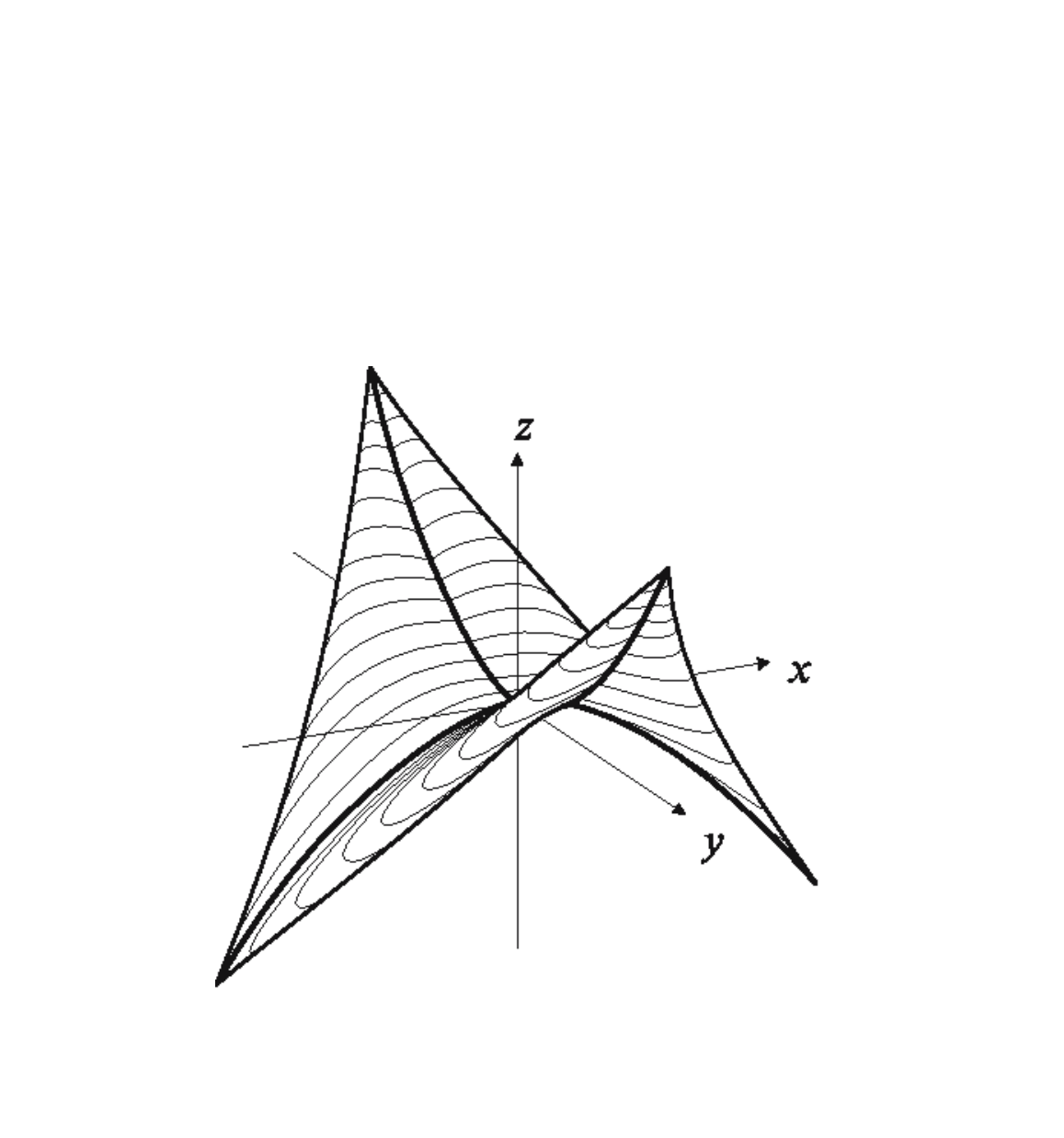}
\caption{The hyperbolic hedgehog of Martinez-Maure}
\label{Hedgehog}
\end{center}
\end{figure}

A second obstruction is that there are no nontrivial one-homogeneous solutions to {\it degenerate} linear elliptic equations in $\mathbb{R}^3$, which are analytic away from the origin. By this we mean there are no nonlinear one-homogeneous functions on $\mathbb{R}^3$, analytic on $\mathbb{S}^2$, such that the two eigenvalues $\lambda_1,\,\lambda_2$ of the Hessian pointwise satisfy either $\lambda_1\lambda_2 < 0$ or $\lambda_1=\lambda_2 = 0$. Alexandrov proved this in \cite{A}, and conjectured that the same should hold if ``analytic" is relaxed to ``smooth". This problem remained open for a while, with several incorrect proof attempts, until Martinez-Maure constructed a surprising counterexample in 2001 \cite{MM}. The hyperbolic hedgehog of Martinez-Maure is smooth away from four cusps that are non-coplanar (Figure \ref{Hedgehog}). We conjecture that the support function of this example minimizes a degenerate convex functional in three dimensions, where $K$ consists of the tips of the four cusps. Such a result would address the third point above (that $K$ being finite does not suffice for $C^1$ regularity of minimizers in dimensions $n \geq 3$), and illustrate the sharpness of the known regularity results ($C^1$ regularity when $K$ has three or fewer points).

\subsection{Exercises}
In this challenging exercise, assume that $u$ is one-homogeneous on $\mathbb{R}^3$, smooth away from the origin, and solves a uniformly elliptic equation of the form
$$\text{tr}(a(x)D^2u) = a_{ij}(x)u_{ij} = 0$$
away from the origin.

\vspace{3mm}

\begin{enumerate}[1.]
\item Show that $w(x_1,\,x_2) : = u(x_1,\,x_2,\,1)$ solves $b_{ij}(x_1,\,x_2)w_{ij} = 0,$
where
$$b_{ij} = (s^Tas)_{ij},\, i,\,j \leq 2; \quad s := I - (0,\,0,\,1) \otimes (x_1,\,x_2,\,0).$$
In particular, if the eigenvalues of $a$ are in $[\lambda,\, \lambda^{-1}]$, then the eigenvalues of $b$ are in $[\lambda,\, (1+ x_1^2+x_2^2)\lambda^{-1}]$,
so the equation for $w$ is locally uniformly elliptic. Conclude that if $\partial_iw$ achieves its maximum somewhere in $\mathbb{R}^2$, then $w$ is linear.

\vspace{3mm}

\item Assume that $q \in \nabla u(\mathbb{S}^2)$ satisfies that $q_3 = \max_{\mathbb{S}^2}\partial_3u$. Show that if $u$ is not linear, then $(\nabla u)^{-1}(q)$ is either the north pole or the south pole. Hint: use the previous result, appropriately rotated. Note that $\nabla u$ is constant on radial lines.

\vspace{3mm}

\item Let $N$ denote the north pole, and assume that $\nabla u(N) = 0$ and $D^2u(N) \neq 0$. Show for $\epsilon$ small that $\nabla u(B_{\epsilon}(N))$ is a smooth graph of the form $\{x_3 = H(x_1,\,x_2)\}$ with $H(0) = 0,\, \nabla H(0) = 0$ and $\det D^2H(0) < 0$. Hint: for $w$ as in Problem $1$, show using the one-homogeneity of $u$ that $H(\nabla w(x)) = w - x\cdot \nabla w$, thus $\nabla H(\nabla w(x)) = -x$. Here $x = (x_1,\,x_2)$.

\vspace{3mm}

\item Assume that $u$ is not linear. Then after a rotation and subtracting a linear function you may assume that $\nabla u(N) = 0$ and $D^2u(N) \neq 0$. Using the previous two parts, show that $\nabla u$ maps south pole to at least two points. Conclude from this contradiction that $u$ must be linear.

\vspace{3mm}

\item Show that the function
$$v(x_1,\,x_2,\,x_3,\,x_4) = \frac{x_1^2+x_2^2-x_3^2-x_4^2}{|x|}$$
solves a linear uniformly elliptic equation of the form $a_{ij}(x)v_{ij} = 0$ in $\mathbb{R}^4$.
\end{enumerate}

\newpage
\section{Open Problems}
\begin{enumerate}[1.]
\item Verify that Martinez-Maure's example from \cite{MM} gives rise to a Lipschitz but non-$C^1$ minimizer in dimension $n = 3$ where $K$ consists of four points. Systematic ways of building hyperbolic hedgehogs in three dimensions
have since been developed using ideas from combinatorial geometry (see \cite{P}), and it would be interesting if these could give rise to a systematic way of building counterexamples to regularity.

\vspace{3mm}

\item Construct parabolic versions of the above examples. In a similar vein, determine whether the singularities in the examples disappear when their boundary data or the integrand $F$ is perturbed.

\vspace{3mm}

\item Study functionals with special structure and symmetry. For example, $F(x) = |x|^p$ with $1 < p < \infty$: is there a sharp estimate on the modulus of continuity of the gradient for Lipschitz minimizers? This problem is well-understood in two dimensions (see e.g. \cite{IM}), but widely open in higher dimensions. Another example is $F(x) = \sum_{i = 1}^n |x_i|^{p_i}$, with $p_i \in (1,\,\infty)$. Here $K$ consists of coordinate hyperplanes when $p_i \neq 2$. The $C^1$ regularity of Lipschitz minimizers is true in two dimensions (see \cite{BB}, \cite{B}), but seems to be open in higher dimensions.

\vspace{3mm}

\item Similar themes arise in parametric geometric variational problems. For example, consider functionals of the form $J(\Sigma) = \int_{\Sigma} \Phi(\nu),$ where $\Sigma \subset \mathbb{R}^{n+1}$ is an oriented hypersurface with unit normal $\nu$, and $\Phi$ is one-homogeneous and convex. Regularity questions for critical points of such functionals (along with their higher-codimension analogues) have attracted recent attention (see e.g. \cite{DDG}, \cite{DT}, \cite{M1}), and it would be interesting to investigate applications of the ideas in the non-parametric setting to such questions.
\end{enumerate}

\newpage

\vspace{5mm}

{\bf Conflict of Interest Statement:} On behalf of all authors, the corresponding author states that there is no conflict of interest. 

\vspace{5mm}


\end{document}